\documentclass[11pt]{article}
\usepackage{amsmath,amscd,amsfonts, amssymb, mathrsfs,textcomp,graphicx}
\usepackage[all]{xy}
\usepackage[russian]{babel}
\usepackage{hyperref}
\usepackage{enumerate}
\usepackage{booktabs}
\begin{document}
\title{{\bfseries Derived $(\infty,1)$-categories of two kinds}}
\author{Grigory Kondyrev}
\date{}
\maketitle 
\newcommand{\Ch}{\operatorname{\sf Ch}}
\newcommand{\PreSp}{\operatorname{\sf PreSp}}
\newcommand{\SSet}{\operatorname{\sf SSet}}
\newcommand{\Hom}{\operatorname{\sf Hom}}
\newcommand{\DGFunct}{\operatorname{\sf DGFunct}}
\newcommand{\Funct}{\operatorname{\sf Funct}}
\newcommand{\Hot}{\operatorname{\sf Hot}}
\newcommand{\Ab}{\operatorname{\sf Ab}}
\newcommand{\Top}{\operatorname{\sf Top}}
\newcommand{\pt}{\operatorname{\sf pt}}
\newcommand{\SigmaSp}{\operatorname{\sf \Sigma-\Sp}}
\newcommand{\OmegaSp}{\operatorname{\Omega-\Sp}}
\newcommand{\N}{\operatorname{\sf N}}
\newcommand{\h}{\operatorname{\sf h}}
\newcommand{\I}{\operatorname{\sf I}}
\newcommand{\AG}{\operatorname{\sf AG}}
\newcommand{\DG}{\operatorname{\sf DG}}
\newcommand{\CDG}{\operatorname{\sf CDG}}
\newcommand{\D}{\operatorname{\sf D}}
\newcommand{\coD}{\operatorname{\sf D^{co}}}
\newcommand{\ctrD}{\operatorname{\sf D^{ctr}}}
\newcommand{\Pro}{\operatorname{\sf Proj}}
\newcommand{\DegProj}{\operatorname{\sf DegProj}}
\newcommand{\DegInj}{\operatorname{\sf DegInj}}
\newcommand{\Cone}{\operatorname{\sf Cone}}
\newcommand{\Acycl}{\operatorname{\sf Acycl}}
\newcommand{\coAcycl}{\operatorname{\sf Acycl^{co}}}
\newcommand{\ctrAcycl}{\operatorname{\sf Acycl^{ctr}}}
\newcommand{\HoProj}{\operatorname{\sf HoProj}}
\newcommand{\HoInj}{\operatorname{\sf HoInj}}
\newcommand{\mmod}{\operatorname{\sf mod}}
\newcommand{\comod}{\operatorname{\sf comod}}
\newcommand{\ctrmod}{\operatorname{\sf ctrmod}}
\newcommand{\Mod}{\operatorname{\sf Mod}}
\newcommand{\Cat}{\operatorname{\sf Cat}}
\newcommand{\Ker}{\operatorname{\sf Ker}}
\newcommand{\Coker}{\operatorname{\sf Coker}}
\newcommand{\Sp}{\operatorname{\sf Sp}}
\newcommand{\ConSp}{\operatorname{\sf ConSp}}
\def\blacksquare{\hbox{\vrule width 5pt height 5pt depth 0pt}}
$${\sf Abstract}$$
The aim of this paper is to reformulate the theory of unbounded derived categories, including more recent categories of first and second kind (following {$ \sf [1]$}), using the language of $(\infty,1)$-categories.
\\
$${\sf Contents}$$
In section 1, we give a brief introduction to the theory of $(\infty,1)$-categories, focusing mainly on the ideas needed further in this paper.
\\
In section 2, we recall the theory of $\DG$-categories. We describe various semiorthogonal decompositions arising in the description of coderived and contraderived categories of modules over a ring and comodules/contramodules over a coring. We also formulate a comodule-contramodule correspondence.
\\
In section 3, we provide a method to construct an $(\infty,1)$-category from a $\DG$-category and investigate some of the main properties of the construction. We also prove localization theorems, which allow us to reformulate main theorems from section 2 in a higher categorical language.
\\
In section 4, we briefly describe curved differential graded analogs for the concepts from section 2 which allows us to easily obtain results similar to those from section 2.
\\
\\
\newpage
\setcounter{section}{-1}
\section{Introduction.}
Let $\mathscr{A}$ be a category of complexes of modules over some ring $R$. 
\\
For the purposes of homological algebra and homotopy theory it is then natural to consider its derived category $\D(\mathscr{A})$ which is defined as the quotient $\D(\mathscr{A})=\h\mathscr{A}/\h \Acycl(A)$ where $\Acycl(A)$ is the full subcategory of $\mathscr{A}$ formed by acyclic complexes of modules and $\h(\bullet)$ is the homotopy category of a given category. Under some finiteness conditions (for example, in the bounded case) it so happens that $\D(\mathscr{A})$ has another usefull description: it is equivalent to the homotopy category of the full subcategory of $\mathscr{A}$ formed by the complexes of degreewise injective/degreewise projective modules.
\\
Without finiteness conditions, the equivalence above does not hold, and we therefore can define two derived categories: the derived category of the first kind, which is defined as $\D^{\sf I}(\mathscr{A})=\h\mathscr{A}/\h \Acycl(A)$, and the derived category of the second kind $\D^{\sf II}(\mathscr{A})$, which is defined as the homotopy category of the full subcategory of $\mathscr{A}$ formed by the complexes of degreewise injective/ degreewise projective modules. It so happens that the choice between projective and injective modules is also significant, and therefore the picture splits into two dual versions of derived categories of the second kind, called the coderived category $\D^{\sf co}(\mathscr{A})$ and the contraderived category $\D^{\sf ctr}(\mathscr{A})$.
\\
Consider now a more general case: let $R$ be a $\DG$-ring and $C$ be a $\DG$-coring. Let $\mathscr{A}$ be one of the following: the category of $\DG$-modules over $R$, the category of $\DG$-comodules over $C$, the category of $\DG$-contramodules over $C$. For each choice of $\mathscr{A}$ we get three different derived categories with many interesting properties. All these definitions and relations between them, as well as properties of the resulting categories, are studied in {$ \sf [1]$}.
\\
\\
The other part of the picture is the theory of $(\infty,1)$-categories. In this paper we use a model of $(\infty,1)$-categories based on quasi-categories, developed mostly in {$ \sf [2]$}, {$ \sf [3]$} and {$ \sf [4]$}, where they are defined as simplicial sets satisfying inner Kan condition.  The theory of quasi-categories is a way to generalize ordinary category theory and provides an appropriate setting for homotopy theory, homological algebra and interactions between them. 
\\
\\
The purpose of this paper is to reformulate the theory of derived categories of two kinds from {$ \sf [1]$} using that language.
\\
\\
{\bfseries Acknowledgements.} This paper is the author's undergraduate research project (year 3, Department of Mathematics, HSE Moscow) supervised by L.Positselsky. 
\\
The author is grateful to L. Positselsky and to A.Gorinov for their useful comments.
\\
\\
\section{Some prerequisites from higher category theory.}
In this section we give a brief introduction to the theory of $(\infty,1)$-categories, focusing mainly on the ideas needed further in this paper.
\\
Recall that an $(\infty,1)$-category is a simplicial set $K$ which has the following property: for any $0<i<n$, any map $f:  \Lambda^{n}_{i} \rightarrow K$ admits an extension $ \widetilde{f}:\Delta^{n} \rightarrow K$ where $\Delta^n$ is the standart $n$-simplex and $\Lambda^{n}_{i}$ is its $i^{th}$ horn. 
\\
It is convenient to think of vertices of $K$ as of objects, edges (not necessary nondegenerate) of $K$ as of morphisms, faces of $K$ (not necessary nondegenerate) as of 2-morphisms, etc. One can see that because of the conditions above there are all ``compositions'' of all ``$p$-morphisms'' and all ``$q$-morphisms'' are invertible for $q>1$ in $K$ up to higher homotopies.
\\
In practice, many $(\infty,1)$-categories are obtained from simplicially enriched categories using the following
\\
\\
{\bfseries Definition 1.1.} There exists a coherent geometric realization functor $$| \bullet |: \SSet \longmapsto \Cat_{\Delta}$$ where $\Cat_{\Delta}$ is a category of simplicial categories (categories enriched over the category $\SSet$ of simplicial sets). 
\\
The simplicial category $| \Delta^{n} |$ is defined as follows:
\\
$\sf Ob$$(| \Delta^{n} |)=\sf Ob(\sf [n])$, where $\sf [n]$ is a poset category $(0 \rightarrow 1 \rightarrow ... \rightarrow n)$.
\\
$\Hom_{| \Delta^{n} | }( \sf k,  \sf m)= 
\begin{cases}
\varnothing , \text{ if } \sf m> \sf k
\\ 
\widetilde{\N}(\sf P_{k,m}), \text{ if } \sf k \leq m,
\end{cases}
$
\\
\\
where $\widetilde{\sf N}$ is a nerve functor $\widetilde{\sf N}: \Cat \longmapsto \SSet$ and $\sf P_{k,m}$ is a poset of subsets of $[k,m]$ that contain both $k$ and $m$ with the partial order given by inclusion. The composition between maps in $| \Delta^{n} |$ is defined naturally.
\\
\\
{\bfseries Remark 1.2.} It is easy to notice that the simplicial set $\widetilde{\N}(\sf P_{k,m})$ is homeomorphic to a cube.
\\
\\
By the formal nonsense, as the category $\Cat_{\Delta}$ admits small colimits, the functor  $| \bullet |$ extends uniquely from the data above to a colimit-preserving functor $| \bullet |: \SSet \longmapsto \Cat_{\Delta}$.
\\
\\
{\bfseries Remark 1.3.} Formally, the functor $| \bullet |: \SSet \longmapsto \Cat_{\Delta}$ can be defined as the coend: $| S |=\int_{}^{\sf [n] \in \Delta}\sf [n] \times S_{n}$, where $S \in \SSet$.
\\
\\
{\bfseries Definition 1.4.} The {\bfseries homotopy coherent nerve functor $\N$ }
$$
\N: \Cat_{\Delta} \longmapsto \SSet
$$
 is characterized by the formula:
$$
\Hom_{\SSet}(\Delta^{n}, \N(\mathscr{C}))=\Hom_{\Cat_{\Delta}}(| \Delta^{n} |,\mathscr{C}).
$$
for $\mathscr{C} \in \Cat_{\Delta}$.
\\
\\
{\bfseries Remark 1.5.} There exists a general technique frequently called {\bfseries nerve and realization}, which provides two adjoint functors $\sf Nerve: \mathscr{B} \longmapsto \SSet$ and $| \bullet |: \SSet \longmapsto \mathscr{B}$ whenever there is a functor $ \Delta \longmapsto \mathscr{B}$ and $\mathscr{B}$ is sufficiently good. The situation above is a special case of this when $\mathscr{B}$ equals $\Cat_{\Delta}$.
\\
\\
{\bfseries Remark 1.6.} Let $\mathscr{C}$ be a simplicial category. Then:
\\
1) The 0-simplex of $\N(\mathscr{C})$ is just an object of $\mathscr{C}$.
\\
2) The 1-simplex of $\N(\mathscr{C})$ is a morphism in $\mathscr{C}$.
\\
3) The 2-simplex of $\N(\mathscr{C})$ can be identified with the diagram in $\mathscr{C}$:
$$
\xymatrix{
X \ar[dr]_-{f_{XY}} \ar[rr]^-{f_{XZ}}  & & Z \\
& Y \ar[ur]_-{f_{YZ}}
}
$$
together with a path from $f_{XZ}$ to $f_{YZ} \circ f_{XY}$ in $\Hom_{\mathscr{C}}(X,Z)$.
\\
\\
We therefore see that given  $\mathscr{C} \in \Cat_{\Delta}$ we can build a simplicial set $\N(\mathscr{C})$. It is natural then to ask under what conditions this simplicial set is an $(\infty,1)$-category.
\\
\\
{\bfseries Proposition 1.7.} Let  $\mathscr{C} \in \Cat_{\Delta}$ be such that for every pair of objects $X, Y \in \mathscr{C}$ it holds that $\Hom_{\mathscr{C}}(X,Y)$ is a Kan complex. Then $\N(\mathscr{C})$ is an $(\infty,1)-$category.
\\
\\
{\bfseries Proof:} We must show that if $0< i < n$, then $\N(\mathscr{C})$ has the  extension property with respect to the inclusion $ \Lambda^{n}_{i} \subseteq \Delta^{n}$, or, equally, that  $\mathscr{C}$ has the extension propery with repsect to the inclusion $| \Lambda^{n}_{i} | \subseteq  | \Delta^{n} |$. First, we can notice that the objects of  $| \Lambda^{n}_{i} |$ are the same as the objects of $ | \Delta^{n} |$, that is, the elements of the poset category $\sf [n]$. Futher, for $0 \leq p \leq q \leq n$ the simplicial set $\Hom_{| \Lambda^{n}_{i} |}(p,q)$ coincides with $\Hom_{ | \Delta^{n}  |}(p,q) $ unless $p=0$ and $q=n$. Consequently, the extension problem
$$
\xymatrix{
\Lambda^{n}_{i} \ar[d] \ar[r]^-{F} & \N(\mathscr{C}) \\
\Delta^{n} \ar@{-->}[ur]
}
$$
is equivalent to
$$
\xymatrix{
\Hom_{| \Lambda^{n}_{i} |}(0,n) \ar[d] \ar[r] & \Hom_{\mathscr{C}}(F(0),F(n)) \\
\Hom_{| \Delta^{n} |}(0,n) \ar@{-->}[ur]
}
$$
Since the simplicial set on the right is a Kan complex by assumption, it suffices to verify that the left vertical map is anodyne. This follows by inspection: the simplicial set $\Hom_{| \Delta^{n} |}(0,n)$ can be identified with the cube  $(\Delta^{1})^{\{1,2,..., n-1\}}$, and $\Hom_{| \Lambda^{n}_{i} |}(0,n)$ can be identified with the simplicial subset obtained by removing the interior of the cube together with one of its faces.
\blacksquare
\\
\\
There is also a natural way to obtain an ordinary category from simplicialy enriched category:
\\
\\
{\bfseries Definition 1.8.} Let $\mathscr{C} \in \Cat_{\Delta}$. Define the {\bfseries homotopy category $\h\mathscr{C}$ of} $\mathscr{C}$ as following:
\\
1) The objects of $\h\mathscr{C}$ are the objects of $\mathscr{C}$.
\\
2) For $X,Y \in \mathscr{C}$, we define $\Hom_{\h\mathscr{C}}(X,Y)=\pi_0 \Hom_{\mathscr{C}}(X,Y)$.
\\
3) The composition of morphisms in $\h\mathscr{C}$ is induced from the composition of morphisms in $\mathscr{C}$ by applying the functor $\pi_0$.
\\
\\
This provides us a process to obtain an ordinary category from any simplicial set:
\\
\\
{\bfseries Definition 1.9.} Define the {\bfseries homotopy functor}
$$
\xymatrix{
\SSet \ar[r]^-{\h} & \Cat
}
$$
as the composition 
$$
\xymatrix{
\SSet \ar[r]^-{| \bullet |} & \Cat_{\Delta} \ar[r]^-{\h} & \Cat
}
$$
If $\mathscr{C}$ is an $(\infty,1)$-category, we then call $\h\mathscr{C}$ its {\bfseries homotopy category}.
\\
\\
{\bfseries Remark 1.10.} Let $\mathscr{C}$ be an $(\infty,1)$-category. The definition above provides us the following presentation of $\h\mathscr{C}$ by generators and relations:
\\
1) The objects of $\h\mathscr{C}$ are the vertices of $\mathscr{C}$.
\\
2) For every edge $ \phi: \Delta^1 \rightarrow \mathscr{C}$, there is a corresponding morphisms $\overline{\phi}$ in $\h\mathscr{C}$ from $\phi(0)$ to $\phi(1)$.
\\
3) For each $\tau: \Delta^2 \rightarrow \mathscr{C}$, it holds that $\overline{d_2(\tau)} \circ \overline{d_0(\tau)}=\overline{d_1(\tau)}$.
\\
4) For each vertex $X$ of $\mathscr{C}$, the morphism $s_0 X$ is the identity morphism $\sf Id_{X}$.
\\
\\
{\bfseries Remark 1.11.} One can observe that $\h: \SSet \rightarrow \Cat$ is actually left adjoint to the nerve functor $\widetilde{\sf N}: \Cat \rightarrow \SSet$.
\\
\\
{\bfseries Remark 1.12.} If $\mathscr{C}$ is a simplicial category, it is not always true that $\h\mathscr{C} \simeq \h\N(\mathscr{C})$. However, this is always true in the case when $\Hom_{\mathscr{C}}(X,Y)$ is a Kan complex for every $X,Y \in \mathscr{C}$ (for example, it holds in the case of Proposition 1.7.).
\\
\\
We now wish to define an appropriate notion of a limit/colimit in an $(\infty,1)$-category. Recall that in ordinary category theory a limit of a diagram can be defined as a final object of a category of objects over that diagram. This motivates the following
\\
\\
{\bfseries Definition 1.13.} Let $\mathscr{C}, \mathscr{D}$ be $(\infty,1)$-categories. Define the {\bfseries join} $\mathscr{C} \star \mathscr{D}$ of $\mathscr{C}$ and $\mathscr{D}$ as the join of the underlying simplicial sets. It is not hard to prove (see, for example, 1.2.8.3 in {$ \sf [2]$}) that the simplicial set $\mathscr{C} \star \mathscr{D}$ inherits the conditions of $(\infty,1)$-category.
\\
\\
{\bfseries Definition 1.14.} Let $\mathscr{C}, \mathscr{D}$ be $(\infty,1)$-categories, and let $p:\mathscr{C} \rightarrow \mathscr {D}$ be a functor between them (a map of simplicial sets). Define the $(\infty,1)$-{\bfseries overcategory} $D_{/ p}$ over $p$ as following: $$({\mathscr{D}_{/ p}})_{n} =\Hom_{p}(\Delta ^{n} \star \mathscr{C}, \mathscr{D}),$$ where the subscript on the right hand side indicates that we consider only those morphisms $f: \Delta^n \star \mathscr{C} \rightarrow \mathscr{D}$ such that $f |\mathscr{C}=p$. 
\\
The $(\infty,1)$-{\bfseries undercategory} $D_{p/}$ under $p$ is defined dually.
\\
\\
{\bfseries Remark 1.15.} By formal nonsense for every $\mathscr{K} \in \SSet$ there exists an adjunction $\Hom_{\SSet}(\mathscr{K},\mathscr{D}_{/ p})=\Hom_{p}(\mathscr{K} \star \mathscr{C}, \mathscr{D})$.
\\
\\
{\bfseries Definition 1.16.} Let $\mathscr{C}$ be an $(\infty,1)$-category. An object $X \in \mathscr{C}$ (a vertice of $\mathscr{C}$) is called {\bfseries final}, if it is final in the category $\h\mathscr{C}$. The {\bfseries initial} object of $\mathscr{C}$ (if it exists) is defined dually.
\\
\\
{\bfseries Remark 1.17.} There exists an equivalent definition of a final object of $(\infty,1)$-category $\mathscr{C}$, namely, one can call an object $X \in \mathscr{C}$ (a vertice of $\mathscr{C}$) final if the projection $\mathscr{C}_{/ X} \rightarrow \mathscr{C}$ is a trivial fibration of simplicial sets.
\\
\\
{\bfseries Definition 1.18.} Let $\mathscr{I}, \mathscr{C}$ be $(\infty,1)$-categories, and let $p:\mathscr{I} \rightarrow \mathscr{C}$ be a functor between them (a map of simplicial sets). A {\bfseries homotopy limit} of $F$ is a final object of $\mathscr{C}_{/ p}$. Dually, a {\bfseries homotopy colimit} is an initial object of $\mathscr{C}_{p / }$.
\\
\\
{\bfseries Remark 1.19.} The notion of a homotopy limit/colimit seems to be the right extension of a notion of a limit/colimit from an ordinary category theory to $(\infty,1)$-category theory. Moreover, homotopy limits/colimits have all the properties of ordinary limits/colimits after they are formulated accurately in the settings of $(\infty,1)$-category theory. The reader is referred to Chapter 4 in {$ \sf [2]$} for a systematical study of homotopy limits and colimits and their properties.
\\
\\
{\bfseries Important Convention.} Under a limit (colimit) hereafter we will mean a homotopy limit (homotopy colimit).
\\
\\
Now all the $(\infty,1)$-categories we will futher work with arise from an algebraic data in a good way and will carry special properties, namely, they will be stable. To be more precise, we start with the
\\
\\
{\bfseries Definition 1.20.} Let $\mathscr{C}$ be an $(\infty,1)$-category which admits zero object and $f:X \rightarrow Y$ be a morphism in $\mathscr{C}$. A {\bfseries kernel} of the morphism $f$ is a cartesian square:
$$
\xymatrix{
\Ker(f) \ar[r] \ar[d] & X \ar[d]^-{f} \\
0 \ar[r] & Y }
$$
Dually, a {\bfseries cokernel} of the morphism $f$ is a cocartesian square:
$$
\xymatrix{
X \ar[d] \ar[r]^-{f} & Y \ar[d] \\
0 \ar[r] & \Coker(f) }
$$
\\
{\bfseries Definition 1.21.} A $(\infty,1)$-category $\mathscr{C}$ is called {\bfseries stable} if the following conditions hold:
\\
1) There exists a zero object $0 \in \mathscr{C}$.
\\
2) Every morphism in $\mathscr{C}$ admits a kernel and a cokernel.
\\
3) A square in $\mathscr{C}$ is cartesian whenever it is cocartesian.
\\
\\
{\bfseries Remark 1.22.} Let $\Omega B$ be the loop space of $B$ defined by the following cartesian square:
$$
\xymatrix{
\Omega C \ar[d] \ar[r]  & 0 \ar[d] \\
0 \ar[r] & B
}
$$
Since $\mathscr{C}$ is stable, for every $B$ $\in$ $\h\mathscr{C}$ (the homotopy category of $\mathscr{C}$) we have $ B \simeq \Omega \Sigma B$ (where $\Sigma B$ is defined as the colimit of the dual diagram).
\\
Therefore $\h\mathscr{C}$ is enriched over the category of abelian groups: 
$$
\pi_{0} \Hom_{\mathscr{C}}(A,B)=\pi_0 \Hom_{\mathscr{C}}(A, \Omega^2 \Sigma^2 B)=\pi_2 \Hom_{\mathscr{C}}(A,\Sigma^2 B),
$$
where the basepoint is given by zero map. It is easy to see that the group structure depends functorially on $X,Y$ $\in$ $h\mathscr{C}$. Moreover, one can see that $\h\mathscr{C}$ has actually a structure of triangulated category where the shift functor is given by the suspension. See theorem 3.11 in $\sf [3]$ for the proof.
\\
\\
Most of the time it is much easier to check whether the given $(\infty,1)$-category is stable using the following description:
\\
\\
{\bfseries Proposition 1.23.} Let $\mathscr{C}$ be an $(\infty,1)$-category which admits zero object. Then the following is equivalent:
\\
(1) The functor ${\Omega}:\mathscr{C} {\longmapsto} \mathscr{C}$ is equivalence and $\mathscr{C}$ admits finite limits.
\\
(2) The functor ${\Sigma}:\mathscr{C} {\longmapsto} \mathscr{C}$ is equivalence and $\mathscr{C}$ admits finite colimits.
\\
(3) The category $\mathscr{C}$ is stable.
\\
\\
{\bfseries Proof:} See corollary 8.28 in {$ \sf [3]$}. \blacksquare
\\
\\
{\bfseries Remark 1.24.} There is a variety of examples of stable $(\infty,1)$-categories. Maybe the main example is the category of pointed $\Omega$-spectrum $\sf Sp $. There  as well exists a procedure called {\bfseries stabilization}, which can be applied to any $(\infty,1)$-category which admits zero object and finite limits and provide somewhat the closest stable $(\infty,1)$-category to the initial one. Namely, if the $(\infty,1)$-category $\mathscr{C}$ is given, we may define a stable $(\infty,1)$-category $\sf Stab(\mathscr{C})$, obtained as a limit of the diagram
\\ 
$$ 
...\overset{\Omega}{\longmapsto} \mathscr{C} \overset{\Omega}{\longmapsto}  \mathscr{C} \overset{\Omega}{\longmapsto}  \mathscr{C}  \overset{\Omega}{\longmapsto} \mathscr{C}.
$$
\\
taken in $(\infty,1)-\Cat$ (the $(\infty, 1)$-category of $(\infty, 1)$-categories). The fact that the category $\sf Stab(\mathscr{C})$ is stable easily follows from Proposition 1.23. See {$ \sf [3]$} for more details.
\\
\\
We now turn our attention to the notion of localization in ($\infty,1$)-categories.
\\
\\
{\bfseries Definition 1.25.} Let $\mathscr{C}$, $\mathscr{D}$ be ($\infty,1$)-categories and $\xymatrix{F,G: \mathscr{C} \ar[r] & \mathscr{D}}$ be functors between them. A {\bfseries natural transformation} from $F$ to $G$ is a simplicial homotopy:
$$\xymatrix{
\mathscr{C} \ar[d]_-{i_0} \ar[rd]^-{F} \\
\mathscr{C} \times \Delta^1 \ar[r] & \mathscr{D} \\
\mathscr{C} \ar[u]^-{i_1} \ar[ur]_-{G}
}$$
\\
\\
{\bfseries Remark 1.26.} One may also proof that a simplicial set of maps between ${\sf Maps}_{\SSet}(\mathscr{C},\mathscr{D})$ between ($\infty,1$)-categories $\mathscr{C}$ and $\mathscr{D}$ is also an ($\infty,1$)-category and then define the natural transformation between the functors $F, G \in {\sf Maps}_{\SSet}(\mathscr{A},\mathscr{D})$ simply as a morphism between them in this category. The reader is referred to Proposition 1.2.7.3. in {$ \sf [2]$} for a more accurate description.
\\
\\
{\bfseries Definition 1.27.} Let $\mathscr{C}$, $\mathscr{D}$ be ($\infty,1$)-categories and $$\xymatrix{\mathscr{C} \ar@/^/[rr]^-{F} && \ar@/^/[ll]^-{G} \mathscr{D}}$$ be functors between them. We say that $F$ is {\bfseries left adjoin} to $G$ if there exists a natrural transformation between $\xymatrix{\epsilon: {\sf Id_{\mathscr{C}}} \ar[r] & G \circ F}$ such that for every $c \in \mathscr{C}$ and $d \in \mathscr{C}$ the map
$$\xymatrix{
\Hom_{\mathscr{D}}(F(c),d) \ar[r] & \Hom_{\mathscr{C}}( G \circ F (c), G(d) ) \ar[r]^-{\epsilon^*} & \Hom_{\mathscr{C}}(c,G(d))
}$$
is an equivalence.
\\
\\
{\bfseries Definition 1.28.} Let $\mathscr{C}$ be an $(\infty,1)$-category and $\mathscr{C}_0$ be its full $(\infty,1)$-subcategory. Than a functor $\xymatrix{L: \mathscr{C} \ar[r] & \mathscr{C}_0}$ is called a {\bfseries localization} if it is left adjoint to the inclusion functor $\xymatrix{i: \mathscr{C}_0 \ar[r] & \mathscr{C}}$. A {\bfseries colocalization} is defined dually.
\\
\\
{\bfseries Proposition 1.29.} Let $\mathscr{C}$ be an $(\infty,1)$-category and  $\xymatrix{i: \mathscr{C}_0 \ar[r] & \mathscr{C}}$  be its full $(\infty,1)$-subcategory. Then the following conditions are equivalent:
\\
1) For every object $X \in \mathscr{C}$ there exists an object $Y \in \mathscr{C}_0$ together with a morphism $\xymatrix{X \ar[r] & i(Y)}$ in $\mathscr{C}$ such that the induced composition $$\xymatrix{\Hom_{\mathscr{C}}(i(Y),i(E)) \ar[r] & \Hom_{\mathscr{C}}(X,i(E))}$$ is a homotopy equivalence for every $E \in \mathscr{C}_0$.
\\
2) The inclusion  $i$ admits a left adjoint.
\\
\\
{\bfseries Proof:} Proposition 5.2.7.8. in {$ \sf [2]$}. \blacksquare
\\
\\
The last thing we will futher need is the following
\\
\\
{\bfseries Definition 1.30.} Let $\mathscr{C}$, $\mathscr{D}$ be ($\infty,1$)-category and $\xymatrix{F: \mathscr{C} \ar[r] & \mathscr{D}}$ be a functor between them. Then $F$ is called an {\bfseries equivalence of ($\infty,1$)-categories} if $F$ is:
\\
1) {\bfseries Full and faithfull}, that is, the map
$$
\xymatrix{F_{x,y}: \Hom_{\mathscr{C}}(x,y) \ar[r] & \Hom_{\mathscr{D}}(F(x),F(y))}
$$
is equivalence of simplicial sets.
\\
2) {\bfseries Essentially surjective}, that is, an induced functor on homotopy categories
$$
\xymatrix{\h F: \h\mathscr{C} \ar[r] & \h\mathscr{D}}
$$
is an essentially surjective functor on ordinary categories.
\\
\\
{\bfseries Remark.} It is worth to say that the definition above is in truth a criterion when a functor between two $(\infty,1)$-categories is an equivalence, just as in ordinary category theory. The definition itself looks quite different, but we will futher not need it in its original form.

\section{$\DG$-categories.}
In this section we first quickly recall the theory of $\DG$-categories. Then we describe various semiorthogonal decompositions arising in the description of coderived and contraderived categories of modules over a ring and comodules or contramodules over a coring. We also formulate a comodule-contramodule correspondence.
\\
\\
{\bfseries Notation.} Let $\mathscr{A}$ a be $\DG$-category, that is, the category enriched over the category of chain complexes $\Ch(\AG)$, where ${\AG}$ is the category of abelian groups. For $C_{\bullet} \in \Ch(\AG)$ denote as $C[n]$ the complex $C_{\bullet-n}$.
\\
\\
{\bfseries Definition 2.1.} Let $A,B \in \mathscr{A}$. A {\bfseries closed morphism} $f:A \rightarrow B$ in $\mathscr{A}$ is an element $f \in \Hom_{0}(A,B)$ such that $d(f) = 0$. The category whose objects are the objects of $\mathscr{A}$ and whose morphisms are closed morphisms is denoted as ${\sf Z}_{0}(\mathscr{A})$.
\\
\\
{\bfseries Definition 2.2.} For $A$ in $\mathscr{A}$, an object $C$ is called the {\bfseries shift} of $A$ by integer $i$ (notation: $C=A[i]$) if a closed isomorphism of contravariant DG-functors $\Hom_{\mathscr{A}}(-,C) \simeq \Hom_{\mathscr{A}}(-,A)[i]$ is fixed.
\\
\\
{\bfseries Definition 2.3.} Let $A,B \in \mathscr{A}$, and let $f: A \rightarrow B$ be a closed morphism. The object $C$ is called the {\bfseries cone} of $f$ (notation: $C=\Cone(f)$) if a closed isomorphism of contravariant DG-functors $\Hom_{\mathscr{A}}(-, C) \simeq \Cone(f_{*})$ is fixed, where $f_{*} : \Hom_{\mathscr{A}}(-,A) \rightarrow \Hom_{\mathscr{A}}(-, B)$.
\\
\\
{\bfseries Definition 2.4.}  Let $\mathscr{A}$ be a $\DG$-category. Denote as $\h \mathscr{A}$ or {\bfseries H}$_0 (\mathscr{A})$ {\bfseries the homotopy category of} $\mathscr{A}$, that is, the category whose objects are the objects of $\mathscr{A}$, and for $A, B \in \h\mathscr{A}$ the morphisms are characterized by the equation $\Hom_{\h\mathscr{A}}(A,B)=${\bfseries H}$_0(\Hom_{\mathscr{A}}(A,B))$.
\\
\\
{\bfseries Definition 2.6.} A $\DG$-{\bfseries ring} $(A,d)$ is pair consisting of associative graded ring $A=\bigoplus_{i \in \mathbb{Z}} A^i$ and an odd derivation $\xymatrix{d: A \ar[r] & A}$ of degree 1 such that $d^2=0$. A $\DG$-{\bfseries coring} is defined dually.
\\
\\
{\bfseries Definition 2.7.} A left {\bfseries module} $M$ over a $\DG$-ring $R$ is an object $M \in \Ch({\AG})$ together with an {\bfseries action} morphism $\xymatrix{ R \otimes M \ar[r]^-{\mu} & M}$ such that the following diagrams commute:
$$\xymatrix{R \otimes R \otimes M \ar[r]^-{Id_{R} \otimes \mu} \ar[d]^-{\bullet_R \otimes Id_M} & R \otimes M \ar[d]^-{\mu} \\
R \otimes M \ar[r]^-{\mu} & M
}$$
$$\xymatrix{
1 \otimes M \ar[dr]_-{\simeq} \ar[rr]^-{\epsilon \otimes Id_{M}}  & & \ar[dl]^-{\mu} R \otimes M \\
& M
}$$
where $\xymatrix{\bullet_R: R \otimes R \ar[r] & R}$ is a multiplication in $R$.
\\
\\
{\bfseries Definition 2.8.} A left {\bfseries comodule} $N$ over a $\DG$-coring $C$ is an object $N \in \Ch({\AG})$ together with a {\bfseries coaction} morphism $\xymatrix{N \ar[r]^-{\eta} & C \otimes N}$ such that the following diagrams commute:
$$\xymatrix{ N \ar[r]^{\eta} \ar[d]^{\eta} & C \otimes N \ar[d]^{\bullet_C \otimes Id_N} \\
C \otimes N \ar[r] ^-{Id_C \otimes {\eta}} & C \otimes C \otimes N
}$$
$$\xymatrix{
C \otimes N \ar[rr]^-{\delta \otimes Id_N}  & &  1 \otimes N \\
& N \ar[ul]^{\mu} \ar[ur]_-{\simeq}
}$$
where $\xymatrix{\bullet_C: C \ar[r] & C \otimes C}$ is a comultiplication in $C$.
\\
\\
{\bfseries Definition 2.9.} A left {\bfseries contramodule} $P$ over a $\DG$-coring $C$ is an object $P \in \Ch({\AG})$ together with a {\bfseries contraaction} morphism $\xymatrix{\Hom(C,P) \ar[r]^-{\alpha} & P}$ such that the following diagrams commute:
$$\xymatrix{ \ar[d]^{\simeq}\Hom(C ,\Hom(C,P)) \ar[rrr]^-{\Hom(C,\alpha)} &&& \Hom(C,P) \ar[d]^-{\alpha} \\
\Hom(C \otimes C, P) \ar[rr]^-{\Hom(\bullet_C,P)} && \Hom(C,P) \ar[r]^-{\alpha} & P
}$$
$$\xymatrix{
\Hom(1,P) \ar[dr]_-{\simeq}  \ar[rr]^-{\Hom(\delta, P)} && \Hom(C,P) \ar[dl]^-{\alpha} \\
& P
}$$
where $\xymatrix{\bullet_C: C \ar[r] & C \otimes C}$ is a comultiplication in $C$.
\\
\\
{\bfseries Remark 2.10.} If one try to formally dualize the notion of a comodule over a coalgbra he or she will end up with the notion of a module over an algebra. If one try to formally dualize the notion of a contramodule over a coalgebra he or she will as well end up with the notion of a module over an algebra. This assymetry arises from the fact that we have an equivalence $\Hom(A \otimes B, C) \simeq \Hom(B, \Hom(A,C))$ but we do not have the dual one.
\\
\\
{\bfseries Convention.} We will futher assume that $\mathscr{A}$ is either of the following:
\\
$\bullet$ A category $R-\mmod$ of left modules over a $\DG$-ring $R$.
\\
$\bullet$ A category $C-\comod$ of left comodules over a $\DG$-coring $C$.
\\
$\bullet$ A category $C-\ctrmod$ of left contramodules over a $\DG$-coring $C$.
\\
One can observe that in either case we can define the complex of homomorphisms of two objects in our category, and therefore $\mathscr{A}$ is a $\DG$ category. One can also observe that in either of the cases in our $\DG$-categories there are products, cones of morphisms and etc.
\\
\\
{\bfseries Definition 2.11.} An object $X \in \mathscr{A}$ is called {\bfseries acyclic}, if $H_i(X)=0$ for every $i$. The full subcategory of $\mathscr{A}$ spanned by acylic objects will be denoted as $\Acycl(\mathscr{A})$.
\\
\\
{\bfseries Remark 2.12.} Since finite direct sum and a cone of a morpism of acyclic objects are again acyclic, as well as the zero object, the homotopy category of acyclic objects $\h\Acycl(\mathscr{A})$ forms a full triangulated subcategory of $\h(\mathscr{A})$.
\\
\\
Now recall from standart homological algebra that if $R$ is an ordinary ring then the bounded derived category of $R$ is equivalent to the homotopy category of complexes of projective/injective $R$-modules.
\\
In our more general case of category of modules (not necessary bounded) over $\DG$-ring or comodules/contramodules over $\DG$-coring this description fails and it is therefore actually possible to define {\bfseries two} derived categories.
\\
We start investigate the first one, namely, we give the following
\\
\\
{\bfseries Definition 2.13.} Define the {\bfseries derived category} of the category $\mathscr{A}$ as following:
$$
\D(\mathscr{A}):=\h(\mathscr{A}) / \h\Acycl(\mathscr{A}).
$$
{\bfseries Remark 2.14.}  The description of the derived category given above as a localization is sometimes rather inconvenient, and it is therefore will be good to have a more explicit description of $\D(\mathscr{A})$.
\\
\\
In order to give another presentation of $\D(\mathscr{A})$, it is usefull to recall the
\\
\\
{\bfseries Definition 2.15.} A {\bfseries semiorthogonal decomposition} $\mathscr{H}=<\mathscr{C},\mathscr{D}>$ of a  triangulated category $\mathscr{H}$ consists of two full triangulated subcategories $\mathscr{C}$, $\mathscr{D}$ and $\mathscr{H}$, such that:
\\
1) For every $C \in \mathscr{C}$, $D \in \mathscr{D}$ it holds that $\Hom_{\mathscr{H}}(C,D)=0$.
\\
2) For every $H \in \mathscr{H}$ there exists a distinguished triangle $C \rightarrow H \rightarrow D$, where  $C \in \mathscr{C}$ and $D \in \mathscr{D}$.
\\
\\
Semiorthogonality is particulary nice in case when we consider localizations of triangulated category. That is, we have a
\\
\\
{\bfseries Proposition 2.16.} Let $\mathscr{H}$ be a triangulated category, and $\mathscr{H}=<\mathscr{C},\mathscr{D}>$ be its semiorthogonal decomposition. Then: 
\\
1) The functor $\mathscr{C} \rightarrow \mathscr{H} / \mathscr{D}$ is an equivalence of categories.
\\
2) The functor $\mathscr{D} \rightarrow \mathscr{H} / \mathscr{C}$ is an equivalence of categories.
\\
3) $\mathscr{C}$ and $\mathscr{D}$ generate $\mathscr{H}$ as triangulated categories, that is, any object of $\mathscr{H}$ can be obtained from objects of $\mathscr{C}$ and $\mathscr{D}$ by iterating the operations of shifts and cone.
\\
4) $\mathscr{C}$ is a full subcategory of $\mathscr{H}$ formed by all objects $C \in \mathscr{H}$ such that $\Hom_{\mathscr{H}}(C,D)=0$ for every $D \in \mathscr{D}$, and an embedding functor $\mathscr{C} \rightarrow \mathscr{H}$ has a right adjoint functor which can be identified with the localization functor $\mathscr{H} \rightarrow \mathscr{H} / \mathscr{D} \simeq \mathscr{C}$.
\\
4) $\mathscr{D}$ is a full subcategory of $\mathscr{H}$ formed by all objects $D \in \mathscr{H}$ such that $\Hom_{\mathscr{H}}(C,D)=0$ for every $C \in \mathscr{C}$, and an embedding functor $\mathscr{D} \rightarrow \mathscr{H}$ has a left adjoint functor which can be identified with the localization functor $\mathscr{H} \rightarrow \mathscr{H} / \mathscr{C} \simeq  \mathscr{D}$.
\\
\\
{\bfseries Proof :} {$ \sf [7]$}. \blacksquare
\\
\\
{\bfseries Definition 2.17.} An object $X \in \mathscr{A}$ is called {\bfseries homotopy projective} if for every acyclic object $A \in \mathscr{A}$ it holds that $\Hom_{\h(\mathscr{A})}(X,A)=0$. The full subcategory of $\mathscr{A}$ spanned by homotopy projective objects will be denoted as $\HoProj(\mathscr{A})$. The full subcategory of {\bfseries homotopy injective} modules $\HoInj(\mathscr{A})$ is defined dually.
\\
\\
{\bfseries Theorem 2.18.} Let $R$ be a $\DG$-ring. Let $C$ be a $\DG$-coring {\bfseries over some fixed field $k$}. Then there are the following semiorthogonal decompositions:
$$
\h(R-\mmod)=< \h\Acycl(R-\mmod), \h\HoInj(R-\mmod)>
$$
$$
\h(R-\mmod)=< \h\HoProj(R-\mmod), \h\Acycl(R-\mmod)>
$$
$$
\h(C-\comod)=< \h\Acycl(C-\comod), \h\HoInj(C-\comod)>
$$
$$
\h(C-\ctrmod)=< \h\HoProj(C-\ctrmod), \h\Acycl(C-\ctrmod)>
$$
Moreover,
\\
1) $\h\HoProj(R-\mmod)$ coincides with $<R>_{\oplus}$, that is, the minimal triangulated subcategory of $\h(R-\mmod)$ spanned by $R$ and closed under infinite direct sums.
\\ 
2) $\h\HoInj(R-\mmod)$ coincides with $<\Hom_{\mathbb{Z}}(R,\mathbb{Q}/\mathbb{Z})>_{\sqcap}$, that is, the minimal triangulated subcategory of $\h(R-\mmod)$ spanned by $\Hom_{\mathbb{Z}}(R,\mathbb{Q}/\mathbb{Z})$ and closed under infinite products. 
\\
4) $\h\HoProj(C-\ctrmod)$ coincides with $<\Hom_{k}(C,k)>_{\oplus}$, that is, the minimal triangulated subcategory of $\h(C-\ctrmod)$ spanned by $\Hom(C,k)$ and closed under infinite direct sums.
\\
4) Using additional set-theoretical assumptions (to be precise, in the assumption of Vop{\u e}nka's principle) $\h\HoInj(C-\comod)$ coincides with $<C>_{\sqcap}$, that is, the minimal triangulated subcategory of $\h(C-\comod)$ spanned by $C$ and closed under infinite products.
\\
\\
{\bfseries Proof:} This appears as 1.4, 1.5, 2.4 and 5.5 {$ \sf [1]$}. \blacksquare
\\
\\
{\bfseries Corollary:}   Let $R$ be a $\DG$-ring. Let $C$ be a $\DG$-coring over some fixed field $k$. From the Proposition 2.16. we then have the following equivalences: 
$$\D(R-\mmod)  \simeq \h\HoInj(R-\mmod) \simeq \h\HoProj(R-\mmod) \simeq <\Hom_{\mathbb{Z}}(R,\mathbb{Q}/\mathbb{Z})>_{\sqcap} \simeq <R>_{\oplus} $$
$$ \D(C-\ctrmod) \simeq \h\HoProj(C-\ctrmod) \simeq <\Hom_{k}(C,k)>_{\oplus}$$
$$ \D(C-\comod) \simeq \h\HoInj(C-\comod) \overset{v}{\simeq} <C>_{\sqcap}$$
when $\overset{v}{\simeq}$ holds at least in the assumption of Vop{\u e}nka's principle.
\\
\\
Following {$ \sf [1]$} we now want to go the second way, that is, to define more exotic versions of the derived categories, so called derived categories of the second kind. 
\\
\\
{\bfseries Definition 2.19.} Define the category $\coAcycl(\mathscr{A})$ of {\bfseries coacyclic} objects as the minimal full triangulated subcategory of $\mathscr{A}$ which contains totalizations of exact triples and closed under infinite coproducts.
\\
\\
{\bfseries Definition 2.20.} Define the category $\ctrAcycl(\mathscr{A})$ of {\bfseries contraacyclic} objects as the minimal full triangulated subcategory of $\mathscr{A}$ which contains totalizations of exact triples and closed under infinite products.
\\
\\
{\bfseries Remark 2.21.} Observe that the categories $\Acycl(\mathscr{A})$ , $\coAcycl(\mathscr{A})$ and $\ctrAcycl(\mathscr{A})$ can differ even in rather simple cases. For example, let $k$ be a field and $\mathscr{A}$ be a category of $\DG$-modules over $\DG$-ring $R=k[x]/(x^2)$ where ${\sf deg}(x)=0$. One then have the following three objects in $R-\mmod$ :
$$
\xymatrix{(1): ...\ar[r]^-{x} & R \ar[r]^-{x} & R \ar[r]^-{x} & R \ar[r]^-{x} & R \ar[r]^-{x} &...}
$$
$$
\xymatrix{(2): ... \ar[r]^-{x} & R \ar[r]^-{x} & R \ar[r] & k \ar[r] & 0 \ar[r] & ...}
$$
$$
\xymatrix{(3): ... \ar[r] & 0 \ar[r] & k \ar[r] & R \ar[r]^-{x} & R \ar[r]^{x} & ...}
$$
that is (2) and (3) are the truncations of (1) above and below respectively.
One can prove that:
 $$(1) \in \Acycl(\mathscr{A}), (1) \notin \coAcycl(\mathscr{A}), (1) \notin \ctrAcycl(\mathscr{A});$$
 $$(2) \in \Acycl(\mathscr{A}), (2) \notin \coAcycl(\mathscr{A}), (2) \in \ctrAcycl(\mathscr{A});$$
 $$(3) \in \Acycl(\mathscr{A}), (3) \in \coAcycl(\mathscr{A}), (3) \notin \ctrAcycl(\mathscr{A});$$
\\
{\bfseries Definition 2.22.} Define the {\bfseries coderived category} of $\mathscr{A}$:
$$ \coD(\mathscr{A}):=\h(\mathscr{A}) / \h\coAcycl(\mathscr{A})$$
and the {\bfseries contraderived category} of $\mathscr{A}$:
$$ \ctrD(\mathscr{A}):=\h(\mathscr{A}) / \h\ctrAcycl(\mathscr{A}).$$
\\
{\bfseries Remark 2.23.} Let $R$ be a $\DG$-algebra and $C$ be a $\DG$-coalgebra. We then have defined the following six categories: $\coD(R-\mmod)$, $\ctrD(R-\mmod)$, $\coD(C-\comod)$, $\ctrD(C-\comod)$, $\coD(C-\ctrmod)$ $\ctrD(C-\ctrmod)$. We suggest that not all of them are sensible enough: namely, the categories $\ctrD(C-\comod)$ and $\coD(C-\ctrmod)$ are ''sufficiently close'' to the zero category (see {$ \sf [1]$) and we will therefore further not focus our attention on them.
\\
\\
{\bfseries Remark 2.24.} Observe from Remark 2.21. that the categories defined above can differ. 
\\
\\
Now to see the parallel with the standart case of complexes of injective/projective modules over an ordinary ring $R$, we introduce the 
\\
\\
{\bfseries Definition 2.25.} An object $X \in \mathscr{A}$ is called {\bfseries degree projective} if having forgotten the differential it becomes a projective module/comodule/contramodule over the base ring/coring with forgotten differential. Denote as $\DegProj(\mathscr{A})$ the full subcategory of $\mathscr{A}$ spanned by degreewise projective objects. The full subcategory $\DegInj(\mathscr{A})$ of {\bfseries degree injective} objects is defined dually.
\\
\\
Let $R$ be a $\DG$-ring. Denote as $R^{\#}$ the graded ring $R$ with forgotten differential. We next define the following two condition on $R$:
\\
(*) A countable direct sum of left graded injective $R^{\#}$-modules has finite injective dimension.
\\
(**) A countable direct product of left graded projective $R^{\#}$-modules has finite projective dimension.
\\
\\
{\bfseries Theorem 2.26.} Let $R$ be a $\DG$-ring. Let $C$ be a $\DG$-coring over some fixed field $k$. Then there are the following semiorthogonal decompositions:
$$
\h(R-\mmod) \overset{*}{=}<\h\coAcycl(R-\mmod)), \h\DegInj(R-\mmod)>
$$
$$
\h(R-\mmod) \overset{**}{=}<\h\DegProj(R-\mmod), \h\ctrAcycl(R-\mmod)>
$$
$$
\h(C-\comod)=<\h\coAcycl(C-\comod)), \h\DegInj(C-\comod)>
$$
$$
\h(C-\ctrmod)=<\h\DegProj(C-\ctrmod), \h\ctrAcycl(C-\ctrmod)>
$$
Where $\overset{*}{=}$ holds when the condition (*) holds for $R$ and $\overset{**}{=}$ holds when the condition (**) holds for $R$.
\\
\\
{\bfseries Proof:} 3.7, 3.8 and 4.4 in {$ \sf [1]$}. \blacksquare
\\
\\
{\bfseries Corollary:}   Let $R$ be a $\DG$-ring. Let $C$ be a $\DG$-coring over some fixed field $k$. From the Proposition 2.26. we then have the following equivalences: 
$$\D^{\sf co}(R-\mmod)  \overset{*}{\simeq} \h\DegInj(R-\mmod) $$
$$\D^{\sf co}(R-\mmod)  \overset{**}{\simeq} \h\DegProj(R-\mmod) $$
$$ \D^{\sf co}(C-\comod) \simeq \h\DegInj(C-\comod) $$
$$ \D^{\sf ctr} (C-\ctrmod) \simeq \h\DegProj(C-\ctrmod) $$
Where $\overset{*}{\simeq}$ holds when the condition (*) holds for $R$ and $\overset{**}{\simeq}$ holds when the condition (**) holds for $R$.
\\
\\
We are now aiming towards the {\bfseries comodule-contramodule correspondence}. To prepare we need the following definition:
\\
\\
{\bfseries Definition 2.27.} Let $C$ be a $\DG$-coring over some fixed field $k$, $N$ be a right $C$-comodule and $P$ a left $C$-contramodule. The {\bfseries contratensor product} $N \odot_{C} P$ of $N$ and $P$ over $C$ is a graded vector space, defined as the coequalizer of the two maps $\xymatrix{ N \otimes_{k} \Hom_{k}(C,P) \ar[r] & N \otimes_{k} P}$, one of which is induced by the left contraaction map, while the other one is obtained as the composition of the coaction map on $C$ and evaluation map $\xymatrix{C \otimes_{k} \Hom_{k}(C,P) \ar[r] & P}$.
\\
\\
Let $C$ be a $\DG$-coring over some fixed field $k$.
\\
For $N \in C-\comod$ one can form the space of comodule homomorphisms $\Hom_{C}(C,M)$ which can be enriched with a structure of $C$-contramodule from the right $C$-comodule structure on $C$. In such a way we obtain a functor $\xymatrix{ \Psi_{C}: C-\comod \ar[r] & C-\ctrmod .}$
\\
For $P \in C-\ctrmod$ one can form a contratensor product $C \odot_{C} P$ which can be enriched with a structure of comodule from the left comodule structure on $C \otimes_{k} P$.  In such a way we obtain a functor $\xymatrix{ \Phi_{C}: C-\ctrmod \ar[r] & C-\comod . }$
\\
\\
We then have the following composition functors:
$$\xymatrix{
\h\DegInj(C-\comod) \ar[r]^-{\Psi_{C}} & \h(C-\ctrmod) \ar[r] & \ctrD(C-\ctrmod)
}$$
$$\xymatrix{
\h\DegProj(C-\ctrmod) \ar[r]^-{\Phi_{C}} & \h(C-\comod) \ar[r] & \coD(C-\comod)
}$$
Abusing the notation we receive the following functors between coderived and contraderived categories of $C$:
$$\xymatrix{
\coD(C-\comod) \ar@/^/[rr]^-{\Psi_C} && \ctrD(C-\ctrmod) \ar@/^/[ll]^-{\Phi_{C}}
}$$
\\
{\bfseries Theorem 2.28.} The functors ${\Phi_C}$ and ${\Psi_C}$ are mutually inverse equivalences between the coderived category $\coD(C-\comod)$ and the contraderived category $\ctrD(C-\ctrmod)$.
\\
\\
{\bfseries Proof:} This appears as $5.2$ in {$ \sf [1]$}. \blacksquare
\\
\\
{\bfseries Remark 2.29.} Observe that in {$ \sf [1]$} all the theorems are formulated in a more general case, that is, not over a $\DG$-ring/coring, but over a ${\sf CDG}$-ring/coring (the notation {\sf CDG} here means curved differential graded).

\section{From $\DG$-category to higher category.}
In this section we recall a method to construct an $(\infty,1)$-category from a $\DG$-category and investigate some of the main properties of the construction. We also prove localization theorems, which allow us to reformulate main theorems from section 2 in a higher categorical language.
\\
\\
{\bfseries Construction 3.1.} Let $\mathscr{A}$ be a $\DG$-category. We can then apply the composition of the Dold-Kan correspondence and the forgetful functor:
\\
 $$
\Ch(\Ab)_{\geq 0}  \longmapsto \sf SAG  \longmapsto \SSet
$$
\\
 to the truncated $\Hom (A,B)_{\bullet}$ complex to receive a simplicial set. Under the truncation functor we mean here the right adjoint to the inclusion 
$\Ch(\AG)_{\geq 0} \subset \Ch(\AG)$. Namely,
 $$\tau_{\geq 0}(...\rightarrow C_{1} \rightarrow C_{0} \rightarrow C_{-1} \rightarrow...)= (...\rightarrow C_{1} \rightarrow {\sf Ker}( C_{0} \rightarrow C_{-1}) \rightarrow 0  \rightarrow...).$$
Therefore, $\mathscr{A}$ has a natural structure of simplicial category. 
\\
\\
{\bfseries Remark 3.2.} For $A, B \in \mathscr{A}$, the homotopy group $\pi_{n}(\Hom(A,B))$ can be identified via the Dold-Kan correspondence with the group of chain-homotopy classes of maps from $A$ to $B[n]$.
\\
\\
{\bfseries Proposition 3.3.} $\N(\mathscr{A})$ is an $(\infty,1)$- category.
\\
\\
{\bfseries Proof:} Every simplicial abelian group is a Kan complex, therefore Proposition 1.1. gives the result.\blacksquare
\\
\\
{\bfseries Definition 3.4.} A $\DG$-category $\mathscr{A}$ is called {\bfseries strongly pretriangulated} if it admits a zero object, all shifts of all objects, and all cones of all morphisms.
\\
\\
{\bfseries Remark 3.5.} For $\mathscr{A}$ any $\DG$-category, there exists a strongly pretriangulated $\DG$-category {\sf PreTr}($\mathscr{A}$) and a fully faithfull $\DG$-functor $\mathscr{A} \rightarrow {\sf PreTr}(\mathscr{A})$, such that for any strongly pretriangulated $\DG$-category $\mathscr{B}$ and any $\DG$-functor $F:\mathscr{A} \rightarrow \mathscr{B}$ there exists a unique lift $F': {\sf PreTr}(\mathscr{A}) \rightarrow \mathscr{B}$. See {$ \sf [6]$} for a proof. 
\\
\\
{\bfseries Important Convention.} Under a category we will hereafter mean either an ordinary category, or $(\infty,1)$-category, depending on the context.
\\
\\
{\bfseries Proposition 3.6.} Let $\mathscr{A}$ be a strongly pretriangulated $\DG$-category. Then $\N(\mathscr{A})$ is a stable category. Moreover, the suspension functor is given by shifting by 1.
\\
\\
{\bfseries Proof:} The idea is to check the second condition in Proposition 1.23.
\\
We remind that in our convention for $C_{\bullet} \in \Ch(\AG)$ we denote as $C[n]$ the complex $C_{\bullet-n}$.
\\
 It is obvious that $\N(\mathscr{A})$ has a zero object (since $\mathscr{A}$ has a zero object). We next prove that cokernels exist in $\N(\mathscr{A})$. Let $g: B \rightarrow A$ be a morphism in $\mathscr{A}$. We assert that the following square is cocartesian in $\N(\mathscr{A})$:
$$
\xymatrix{
A \ar[d] \ar[r]^-{\widetilde{g}} & \Cone(g) \ar[d] \\
0 \ar[r] & B[1]
}
$$
Indeed, it suffices to show that for every $T$ the associated diagram of simplicial sets
$$
\xymatrix{
\Hom_{\mathscr{A}}(B[1],T) \ar[d] \ar[r] & \Hom_{\mathscr{A}}(0,T) \ar[d] \\
\Hom_{\mathscr{A}}(\Cone(g),T) \ar[r]^-{\widetilde{g}^{*}} & \Hom_{\mathscr{A}}(A,T) 
}
$$
is cartesian. As the diagram is obviously a pullback (ordinary pullback, not the homotopy), it suffices to prove that $\widetilde{g}^{*}$ is a Kan fibration. This follows from the fact that ${\widetilde{g}^{*}}$ is the map of simplicial sets associated (under the Dold-Kan correspondence) to a map between complexes of abelian groups which is surjective in positive degrees (every surjection of simplicial abelian groups is a Kan fibration).
\\
It is easy now to describe the suspension functor ${\Sigma}:\N(\mathscr{A}) \rightarrow \N(\mathscr{A})$. Observe that for every $A$ in $\N(\mathscr{A})$ we have the following pushout diagram:
$$
\xymatrix{
A \ar[d] \ar[r] & {\sf Cone}({\sf Id_{A}}) \ar[d] \\
0 \ar[r] & A[1]
}
$$
Now as it holds that ${\sf Cone}({\sf Id_{A}}) \simeq 0$, it follows that the suspension functor ${\Sigma}$$:\N(\mathscr{A}) \rightarrow \N(\mathscr{A})$ can be identified with the shift functor $A \rightarrow A[1]$. In particular, we conclude that ${\Sigma}$ is an equivalence of $(\infty,1)$-categories. 
\\
Observe now that $\N(\mathscr{A})$ has cokernels, as {\sf Coker}($f: A \rightarrow B$) exists via the following cocartesian diagram in $\N(\mathscr{A})$:
$$
\xymatrix{
A \ar[d] \ar[r] & \Cone(\Cone(f)[-1] \rightarrow A) \simeq B\ar[d] \\
0 \ar[r] & \Cone(f)
}
$$
and $A \vee B$ as an object of $\h\N(\mathscr{A})$ can be identified with the colimit of the following diagram in $\N(\mathscr{A})$:
$$
\xymatrix{
A[-1] \ar[d] \ar[r]^-{0} &  B\\
0
}
$$
Therefore, using the Proposition 1.23.,  we conclude that $\N(\mathscr{A})$ is stable. \blacksquare
\\
\\
{\bfseries Remark 3.7.} Let $\mathscr{A}$ be a strongly pretriangulated $\DG$-category, and let  $\mathscr{B}$ be a full  strongly pretriangulated $\DG$-subcategory of $\mathscr{A}$. Then observe that the proof of Proposition 3.6. shows that $\N(\mathscr{B})$ is a stable subcategory of $\N(\mathscr{A})$.
\\
\\
There is also a more straight and algebraic way to obtain an $(\infty,1)$-category from $\DG$-category, namely, we have the following
\\
\\
{\bfseries Construction 3.8.} Let $\mathscr{A}$ be a $\DG$-category. Define a {\bfseries differential graded nerve} $\N_{\DG}(\mathscr{A})$ of $\mathscr{A}$ as followng: for each $n \geq 0$ define $\N_{\DG}(\mathscr{A})_n$ to be a set of all ordered pairs $(\{A_i\}_{0 \leq i \leq n}, f_{I})$, where:
\\
1) For $0 \leq i \leq n$, $A_i$ is an object of $\mathscr{A}$.
\\
2) For every subset $I=\{i_{-} < i_{m} < i_{m-1} < ... < i_{1} < i_{+} \}$ of $\{0<...<n\}$, $f_{I}$ is an element of the abelian group $\Hom_{\mathscr{A}}(A_{i_{-}}, A_{i_{+}})_{m}$, satisfying the equation
$$
d(f_{I})=\sum\limits_{1 \leq j \leq m} (-1)^{j+1}(f_{I-\{i_j\}}-f_{\{i_j < ... <i_m < i_+\}} \circ f_{\{i_- < i_1 <... <i_j\}}).
$$
\\
The simplicial maps are also needed to be carefully written, see 1.3.1.6. in {$ \sf [4]$} for this.
\\
\\
{\bfseries Remark 3.9.}  Let $\mathscr{A}$ be a $\DG$-category. Then:
\\
0) A 0-simplex of $\N_{\DG}(\mathscr{A})$ is simply an object of $\mathscr{A}$.
\\
1) A 1-simplex of $\N_{\DG}(\mathscr{A})$ is a closed morphism in $\mathscr{A}$: that is, a pair of objects $A_0, A_1 \in \mathscr{A}$ together with an element $f_{\{0<1\}} \in \Hom_{\mathscr{A}}(A_0,A_1)_0$, satisfying $d(f_{\{0<1\}})=0$.
\\
2) A 2-simplex of $\N_{\DG}(\mathscr{A})$ consists of a triple of objects $A_0, A_1, A_2 \in \mathscr{A}$, a triple of morphisms: 
$f_{\{0<1\}} \in \Hom_{\mathscr{A}}(A_0,A_1)_0$, $f_{\{1<2\}} \in \Hom_{\mathscr{A}}(A_1,A_2)_0$, $f_{\{0<2\}} \in \Hom_{\mathscr{A}}(A_0,A_2)_0$, such that $d(f_{\{0<1\}})=d(f_{\{1<2\}})=d(f_{\{0<2\}})=0$, together with a morphism $f_{0<1<2} \in \Hom(A_0,A_2)_{1}$ satisfying the equation  $d(f_{\{0<1<2\}})=f_{\{0<2\}} - f_{\{1<2\}} \circ f_{\{0<1\}}$.
\\
\\
{\bfseries Proposition 3.10.}  Let $\mathscr{A}$ be a $\DG$-category. Then:
\\
1) $\N_{\DG}(\mathscr{A})$ is an $(\infty,1)$-category.
\\
2) There is a canonical isomorphism $\h\N_{\DG}(\mathscr{A}) \simeq \h\mathscr{A}$.
\\
3) There exists an equivalence of categories $\theta: \N(\mathscr{A}) \rightarrow \N_{\DG}(\mathscr{A})$ (one should remember that on the left side $\mathscr{A}$ is viewed as the simplicial category via the Dold-Kan correspondence).
\\
\\
{\bfseries Proof:} See $1.3.1$ in {$ \sf [4]$} for the proof. \blacksquare
\\
\\
We are now able to restate the results of section $2$ in a higher categorical language.
\\
\\
{\bfseries Definition 3.11.} Let $R$ be a $\DG$-ring and $C$ be a $\DG$-coring. Define the following $(\infty,1)$-derived categories:
$$ \D(R-\mmod):= \N_{\DG}(\HoInj(R-\mmod)) \simeq \N_{\DG}(\HoProj(R-\mmod))$$
$$ \coD(R-\mmod):=  \N_{\DG}(\DegInj(R-\mmod))$$
$$ \ctrD(R-\mmod):= \N_{\DG}(\DegProj(R-\mmod)) $$
$$ \D(C-\comod):= \N_{\DG}(\HoInj(C-\comod))$$
$$ \D(C-\ctrmod):=\N_{\DG}(\HoProj(C-\ctrmod)) $$
$$ \coD(C-\comod):= \N_{\DG}(\DegInj(C-\comod))$$
$$ \ctrD(C-\ctrmod):= \N_{\DG}(\DegProj(C-\ctrmod)) $$
\\
{\bfseries Proposition 3.12.}  Let $\mathscr{A}$ be a $\DG$-category such that there is a semiorthogonal decomposition $\h\mathscr{A}=<\h\mathscr{C},\h\mathscr{D}>$ for some $\DG$-categories $\mathscr{C}$ and $\mathscr{D}$. Then:
\\
1) There exists a localization functor $\xymatrix{L: \N_{\DG}(\mathscr{A}) \ar[r] & \N_{\DG}(\mathscr{D})}$.
\\
2) There exists a colocalization functor $\xymatrix{CL: \N_{\DG}(\mathscr{A}) \ar[r] & \N_{\DG}(\mathscr{C})}$.
\\
\\
{\bfseries Proof:} We will prove the first statement, the second one is proved dually. Let $\xymatrix{i:\N_{\DG}(\mathscr{D}) \ar[r] & \N_{\DG}(\mathscr{A)}}$ be an inclusion of a full ($\infty,1)$-subcategory. By Proposition 1.29. it is sufficient to find for every $A \in  \N_{\DG}(\mathscr{A})$ an object $D \in  \N_{\DG}(\mathscr{D})$ and a map $\xymatrix{f:A \ar[r] & i(D)}$ such that $f$ induces a homotopy equivalence 
$$\xymatrix{
\Hom_ {\N_{\DG}(\mathscr{A})}(i(D),i(E)) \ar[r] & \Hom_ {\N_{\DG}(\mathscr{A})}(A,i(E))
}$$
for every $E \in \N_{\DG}(\mathscr{D})$. From the definition of a semiorthogonal decomposition for every $A \in \mathscr{A}$ there exist $C' \in \mathscr{C}$ and $D' \in \mathscr{D}$ together with a sequence $\xymatrix{C' \ar[r] & A \ar[r]^{f'} & D'}$ such that it forms a distinguished triangle in $\h(\mathscr{A})$ (and therefore a homotopy colimit diagram in $\N_{\DG}(\mathscr{A})$ as it follows from Proposition 3.10. that $\h\N_{\DG}(\mathscr{A}) \simeq \h\mathscr{A}$). Now one can take $D=D'$ and $f=f'$ and the condition above follows from the long cofiber sequence. \blacksquare
\\
\\
{\bfseries Colloraly:} Let $R$ be a $\DG$-ring. Let $C$ be a $\DG$-coring over some fixed field $k$. 
\\
There are the following localization functors:
$$
\xymatrix{\N_{\DG}(R-\mmod) \ar[r]^-{L} & \N_{\DG}(\HoInj(R-\mmod)) =  \D(R-\mmod)}
$$
$$
\xymatrix{\N_{\DG}(C-\comod) \ar[r]^-{L} & \N_{\DG}(\HoInj(C-\comod)) = \D(C-\comod)}
$$
$$
(*) \xymatrix{\N_{\DG}(R-\mmod) \ar[r]^-{L} & \N_{\DG}(\DegInj(R-\mmod)) = \coD(R-\mmod)}
$$
$$
\xymatrix{\N_{\DG}(C-\comod) \ar[r]^-{L} & \N_{\DG}(\DegInj(C-\comod)) =  \coD(C-\comod)}
$$

There are the following colocalization functors:
$$
\xymatrix{\N_{\DG}(R-\mmod) \ar[r]^-{CL} & \N_{\DG}(\HoProj(R-\mmod))  = \D(R-\mmod)}
$$
$$
\xymatrix{\N_{\DG}(C-\ctrmod) \ar[r]^-{CL} & \N_{\DG}(\HoProj(C-\ctrmod)) = \D(C-\ctrmod)}
$$
$$
(**) \xymatrix{\N_{\DG}(R-\mmod) \ar[r]^-{CL} & \N_{\DG}(\DegProj(R-\mmod)) = \ctrD(R-\mmod)}
$$
$$
\xymatrix{\N_{\DG}(C-\ctrmod) \ar[r]^-{CL} & \N_{\DG}(\DegProj(C-\ctrmod)) =  \ctrD(C-\ctrmod)} 
$$
Where $(*)$ holds when the condition (*) holds for $R$ and $(**)$ holds when the condition (**) holds for $R$. \blacksquare
\\
\\
{\bfseries Colloraly.} Let $C$ be a $\DG$-coring over some fixed field $k$. Following Theorem 2.28. we also have an equivalence of $(\infty,1)$-categories:
$$ \coD(C-\comod) \simeq \ctrD(C-\ctrmod)$$

\section{CDG-theory.}
An abbreviation {\bfseries CDG} states for the notation "curved differential graded". 
\\
One can think of $\CDG$-algebras and $\CDG$-coalgebras as of an extension of $\DG$-algebras and $\DG$-coalgebras where we do not demand the squared differential to be zero. $\CDG$-rings and corings appear in deformations of $\DG$-rings/corings, Koszul duality, non-commutative geometry, TQFTs and etc. In this section we briefly define those objects and the corresponding categories of modules/comodules/contramodules. For the full discussion the reader is referred to {$ \sf [1]$}. To give a taste, however, we make the following
\\ 
\\
{\bfseries Definition 4.1.} A {\bfseries CDG-ring} $B=(b,d,h)$ is a triple consisting of:
\\
1) An assositive graded ring $B= \underset{i \in \mathbb{Z}}{\bigotimes} B^i$.
\\
2) A degree $1$ derivation $\xymatrix{d: B \ar[r] & B}$.
\\
3) An "curvature" element $h \in B^{2}$ satisfying the equations $d^2(x)=[h,x]$ for every $x \in B$ and $d(h)=0$.
\\
\\
{\bfseries Definition 4.2.} A left {\bfseries CDG-module} $M$ over $\CDG$-ring $B$ is a graded lef $B$-module endowed with a degree $1$ derivation $\xymatrix{d_{M}: M \ar[r] & M}$ compatible with the derivation on $B$ and such that $d^2_{M}(m)=hm$ for every $m \in M$.
\\
\\
For any left $\CDG$-modules $M$ and $N$ over $B$, the complex of homomorphisms $\Hom_{B}(M,N)$ from $M$ to $N$ over $B$ is constructed exactly as for $\DG$-modules over a $\DG$-ring. It turns out that these formulas still define a complex in the $\CDG$-module case, as the $h$-related terms cancel each other. 
\\
Therefore left $\CDG$-modules over a given $\CDG$-ring $B$ form $\DG$-category which will further be denoted as $B-\mmod$. All shifts, twists, infinite direct sums, infinite direct products and cones exist in the $\DG$-category of $\CDG$-modules.
\\
\\
For the "dual" notions of {\bfseries CDG-coring}, {\bfseries CDG-comodule/contramodule} over a $\CDG$-coring the reader is referred to the monograph {$ \sf [1]$}.
\\
\\
It is easy to see that in the $\CDG$-case there is no good notion of homology and because of this the derived category is not well defined. However, the full subcategories of coacyclic and contraacyclic objects can be defined just as in the $\DG$-case and therefore the coderived and the contraderived categories are also still defined. One can see that they behave mostly similary as in the $\DG$ case. To be precise, define the full subcategories of degreewise injective and degreewise projective objests just as in the $\DG$-case. We then can state the theorem below which is proved similary as in section 3.
\\
\\
{\bfseries Theorem 4.3.} Let $B$ be a $\CDG$-ring and $E$ be a $\CDG$-coring over some fixed field $k$. Then
1) There are the following semiorthogonal decompositions:
$$
\coD(B-\mmod) \overset{*}{=}<\h\coAcycl(B-\mmod)), \h\DegInj(B-\mmod)>
$$
$$
\ctrD(B-\mmod) \overset{**}{=}<\h\DegProj(B-\mmod), \h\ctrAcycl(B-\mmod)>
$$
$$
\coD(E-\comod)=<\h\coAcycl(E-\comod)), \h\DegInj(E-\comod)>
$$
$$
\ctrD(E-\ctrmod)=<\h\DegProj(E-\ctrmod), \h\ctrAcycl(E-\ctrmod)>
$$
2) There is as well an equivalence of categories:
$$\xymatrix{
\coD(E-\comod) \ar@/^/[rr]^-{\Psi_C} && \ctrD(E-\ctrmod) \ar@/^/[ll]^-{\Phi_{C}}
}$$
\\
{\bfseries Proof:} 3.7, 3.8 ,4.4 and 5.2 in {$ \sf [1]$}. \blacksquare
\\
\\
We can now simply apply the differential graded nerve construction and simulate the proof as for $\DG$ case, receiving the corresponding $(\infty,1)$-categories, localization, colocalization and correspondence theorems for $\CDG$ situation. To be precise, we have the following
\\
\\
{\bfseries Theorem 4.4.} Let $B$ be a $\CDG$-ring and $E$ be a $\CDG$-corind over some fixed field $k$. There are the following localization functors:
$$
(*) \xymatrix{\N_{\DG}(B-\mmod) \ar[r]^-{L} & \N_{\DG}(\DegInj(B-\mmod)) = \coD(B-\mmod)}
$$
$$
\xymatrix{\N_{\DG}(E-\comod) \ar[r]^-{L} & \N_{\DG}(\DegInj(E-\comod)) =  \coD(E-\comod)}
$$
There are also the following colocalization functors:
$$
(**) \xymatrix{\N_{\DG}(B-\mmod) \ar[r]^-{CL} & \N_{\DG}(\DegProj(B-\mmod)) = \ctrD(B-\mmod)}
$$
$$
\xymatrix{\N_{\DG}(E-\ctrmod) \ar[r]^-{CL} & \N_{\DG}(\DegProj(E-\ctrmod)) =  \ctrD(E-\ctrmod)} 
$$
Where $(*)$ holds when the condition (*) holds for $B$ and $(**)$ holds when the condition (**) holds for $B$. \blacksquare

\section{References.}
{$ \sf [1]$} {\em L. Positselski}, Two kinds of derived categories, Koszul duality, and comodule-contramodule correspondence,
\\
 arXiv:0905.2621.
\\
\\
{$ \sf [2]$} {\em J. Lurie}, Higher Topos Theory,
\\
 arXiv:math/0608040.
\\
\\
{$ \sf [3]$}  {\em J. Lurie}, Derived Algebraic Geometry I: Stable $\infty$-Categories, 
\\
arXiv:math/0608228.
\\
\\
{$ \sf [4]$}  {\em J. Lurie}, Higher Algebra, 
\\
http://www.math.harvard.edu/~lurie/papers/HigherAlgebra.pdf.
\\
\\
{$ \sf [5]$} {\em C. Weibel}, An Introduction to Homological Algebra,
\\
 Cambridge University Press, 1995.
\\
\\ 
{$ \sf [6]$} {\em A. Bondal, M. Kapranov}, Enhanced Triangulated Categories, 
\\
Math., USSR Sbornik, Vol. 70 (1991).
\\
\\ 
{$ \sf [7]$} {\em J.-L. Verdier}, Categories derivees, etat 0,
\\
SGA 4 1/2. Lect. Notes Math. 569, p. 262-311, 1977.
\\
\\
{$ \sf [8]$} {\em J. Adams}, Stable homotopy and generalised homology, 
\\
University of Chicago Press, Chicago, 1974.
\\
\\
Department of Mathematics, National Research University Higher School of Economics.
\\
E-mail address: gkondyrev@gmail.com
\end{document}